\documentclass[12pt,reqno]{amsart}
\usepackage{amssymb,amsthm,amsmath,xypic}
\usepackage[all]{xy}
\usepackage[dvips]{color}

\newcommand{\F}{{\mathcal F}}

\newcommand{\G}{{\mathcal G}}
\newcommand{\Hc}{{\mathcal H}}
\newcommand{\A}{\mathcal A}
\newcommand{\V}{\mathcal V}
\newcommand{\LC}{\mathcal L}
\newcommand{\CC}{\mathcal C}
\newcommand{\mD}{\mathcal D}
\newcommand{\cal}{\mathcal}
\newtheorem{theorem}{Theorem}[section]
\newtheorem{corollary}[theorem]{Corollary}
\newtheorem{lemma}[theorem]{Lemma}
\newtheorem{proposition}[theorem]{Proposition}
\newtheorem{proposition-definition}[theorem]{Proposition - Definition}
\newtheorem{definition}[theorem]{Definition}
\newtheorem{remark}[theorem]{Remark}
\newtheorem{example}[theorem]{Example}

\begin{document}

\title[Categories as models on a suitable algebraic theory]{Categories as models \\on a suitable algebraic theory}
\author{Kuerak Chung and Giovanni Marelli}

\begin{abstract}
We explain how categories, and groupoids, can be seen as models for a Lawvere ${\mathfrak Gr}$-theory, where ${\mathfrak Gr}$ is the category of graphs, and show that for Lawvere ${\mathfrak Gr}$-theories finitely presentable models are finitely presentable objects.
\end{abstract}

\keywords{Lawvere theories, finitely presentability, MSC2010 18C10}
\address{KIAS,  Hoegiro 87, Dongdaemun-gu, Seoul 130-722, South Korea,\\
 krchung@kias.re.kr\\
KIAS,  Hoegiro 87, Dongdaemun-gu, Seoul 130-722, South Korea,\\
and Instituto de Matematicas, Universidad de Antioquia, Bloque 4, Calle 70 n.52-21, Medellin,
Colombia\\
marelli@kias.re.kr\\
marelli@matematicas.udea.edu.co}
\maketitle

\section{Introduction}
Lawvere theories were introduced by Bill Lawvere in his doctoral thesis \cite{L} in 1963 as a categorical
formulation
of universal algebra. The correspondence between Lawvere theories and finitary monads on $\mathfrak{Set}$ is one of the deepest relationships in category theory. In \cite{P} Lawvere
theories were generalized to enriched Lawvere theories, substituting $\mathfrak{Set}$
with an arbitrary base category $\cal{V}$ satisfying axioms that make $\cal{V}$
an appropriate base category for enrichment in the sense of \cite{K}, and a correspondence between $\cal{V}$-enriched Lawvere theories and $\cal{V}$-enriched monads on $\cal{V}$ was achieved. A further step was taken in \cite{NP} and \cite{LP} with the notion of Lawvere ${\cal A}$-theories:
first a category $\cal{V}$ in which to enrich and then a base $\cal{V}$-category ${\cal A}$ were chosen. The correspondence above was extended to one between Lawvere ${\cal A}$-theories and finitary $\cal{V}$-enriched monads on the $\cal{V}$-cateogry $\cal{A}$. This allowed to
view as models for Lawvere $\cal{A}$-theories structures for which this interpretation was not possible with $\cal{A}=\cal{V}$.

In this paper we first show, as an application of what explained above, that categories and groupoids can be seen as models for certain Lawvere ${\mathfrak Gr}$-theories, where  $\cal{A}={\mathfrak Gr}$ is the category of graphs and $\cal{V}={\mathfrak Set}$.

Another property of Lawvere theories on $\mathfrak{Set}$ is that a model $M$ for a given theory is finitely presentable exactly when $Mod(M,-):Mod\rightarrow{\mathfrak Set}$ preserves filtered colimits, where $Mod$ denotes the category of models for the given theory. This provides an equivalence between
an extrinsic (the former) and an intrinsic (the latter) characterization of finitely presentability. We
show that this still holds for categories, seen, as said, as models for a Lawvere ${\mathfrak Gr}$-theory, where the fact that $\cal{A}={\mathfrak Gr}$ is decisive. We do not know if this equivalence holds for generic Lawvere $\cal{A}$-theories and at the moment we have not counterexamples.

The paper is organized as follows: in the second chapter we remind the notion of graph and resume their basic properties; in the third we remember Lawvere ${\cal A}$-theories, for a locally finitely presentable $\V$-category $\A$, where $\V$ is a locally finitely presentable symmetric monoidal closed category, and their $\V$-category of models, particularly we show how categories and groupoids can be seen each one as models for a suitable Lawvere ${\mathfrak Gr}$-theory, where ${\mathfrak Gr}$ denotes the category of graphs; finally, in the fourth, we show that finitely presentable categories are just finitely presentable models, establishing an equivalence between an intrinsic and extrinsic characterization.

We would like to thank Bernhard Keller, who gave us a motivation for studying this kind of problems, and for useful discussions. We wish to thank also Ross Street, Stephen Lack and John Power for useful explanations and suggestions.

\section{Graphs}

We introduce here the notion of graph, explaining some of their properties, and the category of graphs and graphs morphisms.

\begin{definition}
\label{graph}
A (directed) graph ${\mathcal G}$ consists of
\begin{enumerate}
\item a class ${\G}_0$, whose elements are called vertices (or 0-cells);
\item for each pair $(A,B)\in\G_0\times \G_0$ a set $\G(A,B)$, whose elements are called the arrows (or 1-cells or edges) from $A$ to $B$.
\end{enumerate}
\end{definition}

Equivalently, we can assign a graph $\G$ by giving a class ${\G}_0$ of vertices and a class ${\G}_1$ of arrows, together with two maps of classes $s,t:{\G}_1\rightarrow {\G}_0$, called source and target, such that the arrows with given source and target form a set.

\begin{definition}
A morphism of graphs $\alpha:\G\rightarrow\Hc$ between two graphs $\G$ and $\Hc$ consists of
\begin{enumerate}
\item a map $\alpha_0:{\G}_0\rightarrow{\Hc}_0$
\item for each $(A,B)\in\G\times\G$ a map $\alpha_{A,B}$
$$\G(A,B)\rightarrow\Hc(\alpha A,\alpha B)$$
\end{enumerate}
\end{definition}

Equivalently, a morphism of graphs $\alpha$ is assigned by giving maps $\alpha_0:{\G}_0\rightarrow{\Hc}_0$ and $\alpha_1:{\G}_1\rightarrow{\Hc}_1$ commuting with $s$ and $t$.

\begin{proposition}
Small graphs and morphisms of graphs form a category, which we denote by ${\mathfrak Gr}$.
\end{proposition}

Another useful characterization of graphs is that of presheaves over a suitable category. Let ${\mathfrak Set}$ be the category of sets and $\mD$ is the subcategory of ${\mathfrak Set}$, whose objects are the sets $\bar{0}:=\{0\}$ and $\bar{1}:=\{0,1\}$, and whose non-trivial morphisms are the obvious inclusions $i_0,i_1:\{0\} \rightarrow \{0,1\}$ to $0, 1$ respectively;
\begin{displaymath}
\xymatrix{
\bullet_{\bar{0}} \ar@/^/[r]|{i_0}
        \ar@/_/[r]|{i_1} &
\bullet_{\bar{1}}
}
\end{displaymath}

\begin{proposition}
\label{presheaves}
${\mathfrak Gr}$ is isomorphic to ${\mathfrak Set}^{\mD^{op}}$.
\end{proposition}
\begin{proof}
Given a graph ${\G}=({\G}_0,{\G}_1,s,t)$ we define a presheaf $\Phi$ on $\mD$ by setting $\Phi(\bar{0})={\G}_0$, $\Phi(\bar{1})={\G}_1$, $\Phi(i_0)=s$, $\Phi(i_1)=t$; conversely, the same definitions assign to a given presheaf $\Phi$ a graph $\G$. Given a morphism $\alpha:\G\rightarrow\Hc$, clearly from the equality above, it defines a morhism between presheaves $\Phi$ and $\Psi$ defined by $\G$ and $\Hc$ respectively, and the converse holds too.
\end{proof}

As examples we compute the graphs associated to the representable functors $h_{\bar{0}}(-)=Hom_\mD(-,\bar{0})$ and $h_{\bar{1}}(-)=Hom_\mD(-,\bar{1})$ in ${\mathfrak Set}^{\mD^{op}}$.

\begin{example}
From the definition of $\mD$, we have that $h_{\bar{0}}(\bar{0})=\{id_{\bar{0}}\}$ and $h_{\bar{0}}(\bar{1})=\varnothing$, so that $h_{\bar{0}}$ is the graph with one vertex and no arrows;
\begin{displaymath}
\xymatrix{
\bullet_{id_{\bar{0}}}
}.
\end{displaymath}

Instead $h_{\bar{1}}(\bar{0})=\{i_0,i_1\}$ and $h_{\bar{1}}(\bar{1})=\{id_{\bar{1}}\}$, so that $h_{\bar{1}}$ is a graph with two vertexes and one arrow $id_{\bar{1}}$ from $i_0$ to $i_1$;
\begin{displaymath}
\xymatrix{
\bullet_{i_0} \ar[r]^{id_{\bar{1}}} & \bullet_{i_1}
}.
\end{displaymath}
\end{example}

\begin{corollary}
${\mathfrak Gr}$ is locally finitely presentable.
\end{corollary}
\begin{proof}
It follows from the fact that ${\mathfrak Gr}$ is a category of presheaves by proposition \ref{presheaves}.
\end{proof}

In particular, ${\mathfrak Gr}$ is complete and cocomplete such that limits and colimits can be computed pointwisely, or, equivalently, according to definition \ref{graph}, cellwisely.

The following proposition establishes a relation between the category ${\mathfrak Cat}$ of small categories and the category $\mathfrak{Gr}$ of graphs:

\begin{proposition}
\label{free}
As a functor between ${\mathfrak Set}$-categories, the forgetful functor ${\cal U}:{\mathfrak Cat}\rightarrow{\mathfrak Gr}$ has a left adjoint ${\cal F}$.
\end{proposition}
\begin{proof}
See \cite{B}.
\end{proof}

\begin{remark}
\label{enrich}
\rm ${\mathfrak Gr}$ is a symmetric monoidal closed category. ${\mathfrak Gr}$ and ${\mathfrak Cat}$ are enriched over ${\mathfrak Gr}$,
however proposition \ref{free} does not extend to ${\mathfrak Gr}$-adjunction.
\end{remark}

\section{Lawvere ${\cal A}$-theories}

As explained in remark \ref{enrich} we will be concerned with Lawvere ${\cal A}$-theories when ${\cal A}=\mathfrak{Gr}$ and ${\cal V}={\mathfrak Set}$, however, following
\cite{NP}, we introduce them in generality. Suppose that ${\cal V}$ is locally finitely presentable as a symmetric monoidal closed category and that ${\cal A}$ is a locally finitely presentable ${\cal V}$-category. Denote by ${\cal A}_{fp}$ a skeleton of the full sub-${\cal V}$-category of ${\cal A}$ given by finitely presentable objects of ${\cal A}$. Let $i:{\cal A}_{fp}\rightarrow{\cal A}$ be the inclusion ${\cal V}$-functor and $\tilde{i}$ the following composition:
\begin{displaymath}
\xymatrix{
{\cal A} \ar[r]^Y & [{\cal A}^{op},{\cal V}] \ar[r]^{[i^{op},{\cal V}]} & [{\cal A}^{op}_{fp},{\cal V}] }
\end{displaymath}
where $Y$ is the enriched Yoneda embedding. As to ${\mathfrak Gr}$, note that finitely presentable objects are just finite graphs; we will denote ${\mathfrak Gr}_{fp}$ simply by ${\mathfrak Gr}_f$.
\begin{definition}
A Lawvere ${\cal A}$-theory is a small ${\cal V}$-category ${\cal L}$ together with an identity-on-objects strict finite ${\cal V}$-limit- preserving ${\cal V}$-functor $J:{\cal A}^{op}_{fp}\rightarrow{\cal L}$.
\end{definition}

\begin{definition}
\label{model}
Given a Lawvere ${\cal A}$-theory $({\cal L},J)$, its ${\cal V}$-category of models is defined by the following pull-back in the ${\cal V}-Cat$ of locally small ${\cal V}$-categories:
\begin{displaymath}
\xymatrix{
Mod({\cal L}) \ar[r]^{P_{\cal L}} \ar[d]_{U_{\cal L}} &
[{\cal L},{\cal V}] \ar[d]^{[J,{\cal V}]} \\
{\cal A} \ar[r]_{\tilde{i}}
& [{\cal A}^{op}_{fp},{\cal V}] }
\end{displaymath}
\end{definition}

We quote the following result from \cite{NP}:

\begin{proposition}
\label{monadic}
$U_{\cal L}$ is finitary monadic, particularly it has a left ${\cal V}$-adjoint $F_{\cal L}$
\end{proposition}

For simplicity, when the theory $\LC$ is fixed, we will use the notation $U$ and $F$ for the forgetful functor and its left adjoint.

As said, we want to show that categories can be seen as models for an ${\cal A}$-Lawvere theory with ${\cal V}={\mathfrak Set}$ and ${\cal A}={\mathfrak Gr}$.

Let $\overrightarrow{0}$ be the following graph which is isomorphic to the graph corresponding to the representable functor $h_{\bar{0}}$ in ${\mathfrak Set}^{\mD^{op}}$
 \begin{displaymath}
\xymatrix{
\overrightarrow{0}:= & \bullet_{a}
}
\end{displaymath}
and $\overrightarrow{1}$ the following graph which is isomorphic to the graph corresponding to the representable functor $h_{\bar{1}}$ in ${\mathfrak Set}^{\mD^{op}}$
\begin{displaymath}
\xymatrix{
 \overrightarrow{1}:= & \bullet_{a} \ar[r] & \bullet_{b}
}.
\end{displaymath}
 By abuse of notations, $s$ and $t$ denote the two morphisms of graphs from $\overrightarrow{0}$ to $\overrightarrow{1}$, mapping the only vertex of $\overrightarrow{0}$ to $a$ and $b$ respectively
\begin{displaymath}
\xymatrix{
\bullet_{a} \ar@/^1pc/[rr]^t \ar@/_/[r]_s & \bullet_{a} \ar[r] & \bullet_{b}
}
\end{displaymath}

Note that the graph $\overrightarrow{2}$, defined as the graph with three vertexes $a$, $b$ and $c$ and two arrows from $a$ to $b$ and from $b$ to $c$
\begin{displaymath}
\xymatrix{
\overrightarrow{2}:= & \bullet_{a} \ar[r] & \bullet_{b} \ar[r] & \bullet_{c}
}
\end{displaymath}
is the push-out of $s$ and $t$ in $\mathfrak{Gr}$
\begin{displaymath}
\xymatrix{
\vec{0} \ar[r]^t \ar[d]_s & \vec{1} \ar[d]^{s'} \\
\vec{1} \ar[r]_{t'} & \vec{2} }
\end{displaymath}
,i.e., $\overrightarrow{2} \cong \overrightarrow{1}+_0\overrightarrow{1}$.
In a similar way, the graph
\begin{displaymath}
\xymatrix{
\overrightarrow{3}:= & \bullet_{a} \ar[r] & \bullet_{b} \ar[r] & \bullet_{c} \ar[r] & \bullet_{d}
}
\end{displaymath}
is isomorphic to  $\overrightarrow{1}+_{\overrightarrow{0}}\overrightarrow{1}+_{\overrightarrow{0}}\overrightarrow{1}$ in $\mathfrak{Gr}$.

In general, $$\overrightarrow{n};=\bullet_{a_0} \rightarrow \bullet_{a_1} \cdots \rightarrow \bullet_{a_n} \cong \overrightarrow{1}+_{\overrightarrow{0}} \cdots +_{\overrightarrow{0}} \overrightarrow{1}$$.

We may consider that above graphs and morphisms are in $\mathfrak{Gr}_f$ and above finite colimits are those in $\mathfrak{Gr}_f$ since $i:\mathfrak{Gr}_f \rightarrow \mathfrak{Gr}$ preserves finite colimits.

Note that for any graph $G$
$$\mathfrak{Gr}(\overrightarrow{0},G) \cong G_0, \mathfrak{Gr}(\overrightarrow{1},G) \cong G_1, \mathfrak{Gr}(\overrightarrow{n},G) \cong G_1 \times_{G_0} G_1 \times_{G_0} \cdots \times_{G_0} G_1$$. In particular, we have the following cartesian (pullback) diagram $\mathfrak{Gr}(\overrightarrow{2},G)$ in $\mathfrak{Set}$ corresponding to the pushout diagram $\overrightarrow{1}+_{\overrightarrow{0}}\overrightarrow{1}$ in $\mathfrak{Gr} $;
\begin{displaymath}
\xymatrix{
G_1 \times_{G_0} G_1 \ar[r]^{t'} \ar[d]_{s'} & G_1 \ar[d]^s\\
G_1 \ar[r]_t  & G_0.}
\end{displaymath}

Denote the obvious inclusions in $\mathfrak{Gr}$ by

$l_j:\overrightarrow{1} \rightarrow \overrightarrow{3}, j=1,2,3$

$l_{jk}:\overrightarrow{2} \rightarrow \overrightarrow{3},(j,k)=(1,2),(2,3)$.

We define now the Lawvere theory we are interested in.
\begin{definition}
\label{cat}
${\LC}_{\mathfrak C}$ is the Lawvere ${\mathfrak Gr}$-theory having the following presentation;

generators:  $ m:\overrightarrow{2}\rightarrow\overrightarrow{1}$, $e:\overrightarrow{0} \rightarrow \overrightarrow{1}$

 axioms(relations):
\begin{displaymath}
\xymatrix{
\vec{2} \ar[r]^{m} \ar[d]_{s'^{op}} & \vec{1} \ar[d]^{s^{op}} & & \vec{2} \ar[r]^{m} \ar[d]_{t'^{op}} & \vec{1} \ar[d]^{t^{op}}& &\vec{3} \ar[r]^{\psi} \ar[d]_{\phi} & \vec{2} \ar[d]^m &\\
\vec{1} \ar[r]_{s^{op}} & \vec{0} &,& \vec{1} \ar[r]_{t^{op}}
& \vec{0} &  &
\vec{2} \ar[r]_m & \vec{1}& , }
\end{displaymath}

\begin{displaymath}
\xymatrix{
\vec{0} \ar[r]^{e} \ar[dr]_{id} & \vec{1} \ar[d]^{s^{op}} &
 &\vec{0} \ar[r]^{e} \ar[dr]_{id} & \vec{1} \ar[d]^{t^{op}} &
 & \vec{1} \ar[r]^{\delta} \ar[dr]_{id} & \vec{2} \ar[d]^m & \vec{1} \ar[l]_{\rho} \ar[dl]^{id}\\
  & \vec{0}  & ,
 & & \vec{0} &,
 &  & \vec{1} &.
 }
\end{displaymath}

where  $\psi,\phi, \delta, \rho$ are the unique morphisms in ${\LC}_{\mathfrak C}$ making the following diagrams in ${\LC}_{\mathfrak C}$ commute
\begin{displaymath}
\xymatrix{
 & & \vec{3} \ar[lld]_{{l_{12}}^{op}} \ar@{.>}[d]|-{\psi} \ar[drr]^{{l_3}^{op}} & &    &   & &
     \vec{3} \ar[lld]_{{l_{1}}^{op}} \ar@{.>}[d]|-{\phi} \ar[drr]^{{l_{23}}^{op}} & & \\
\vec{2} \ar[rd]_m & & \vec{2} \ar[ld]_{{s'}^{op}} \ar[rd]^{{t'}^{op}} & & \vec{1} \ar[dl]^{id}   &
\vec{1} \ar[rd]_{id} & & \vec{2} \ar[ld]_{{s'}^{op}} \ar[rd]^{{t'}^{op}} & & \vec{2} \ar[dl]^{m}  \\
 & \vec{1} \ar[rd]_{t^{op}} & & \vec{1} \ar[dl]^{s^{op}} &  &   &
   \vec{1} \ar[rd]_{t^{op}} & & \vec{1} \ar[dl]^{s^{op}} & \\
 & & \vec{0} & &  ,&  & &
     \vec{0} & &
 }
\end{displaymath}

\begin{displaymath}
\xymatrix{
\vec{1} \cong \vec{0} \times_{\vec{0}} \vec{1} \ar[d]_{s^{op}} \ar@{.>}[dr]|-{\delta} \ar[drr]^{id}   & &
& &\vec{1}  \cong \vec{1} \times_{\vec{0}} \vec{0}\ar[ddr]_{id} \ar[r]^{t^{op}} \ar@{.>}[dr]|-{\rho}& \vec{0} \ar[rd]^e &\\
 \vec{0} \ar[rd]_e & \vec{2} \ar[r]^{{t'}^{op}} \ar[d]_{{s'}^{op}} & \vec{1} \ar[d]^{s^{op}} &
& & \vec{2} \ar[r]^{{t'}^{op}} \ar[d]_{{s'}^{op}} & \vec{1} \ar[d]^{s^{op}} & \\
  & \vec{1} \ar[r]_{t^{op}} & \vec{0}
&, & & \vec{1} \ar[r]_{t^{op}} & \vec{0}
   .}
\end{displaymath}

\end{definition}
Note that such unique morphisms $\psi, \phi, \delta, \rho$ exist in ${\LC}_{\mathfrak C}$, since the bottom diagrams are cartesian in ${\LC}_{\mathfrak C}$ and the outer diagrams commute (by axioms).

The next theorem says that categories are the models for this theory.

\begin{theorem}
\label{catt}
The category $Mod({\LC}_{\mathfrak C})$ of ${\LC}_{\mathfrak C}$-models is equivalent to the category ${\mathfrak Cat}$.
\end{theorem}
\begin{proof}
From definition \ref{model} we have that for any model $M$ there exists a graph $G\in {\mathfrak Gr}$ such that $M\circ J={\mathfrak Gr}(i-,G)$.

The first two diagrams yield the following commutative diagrams in $\mathfrak{Set}$
\begin{displaymath}
\xymatrix{
G_1\times_{G_0} G_1 \ar[r]^{~~~~~~~~~~~~~M(m)} \ar[d]_{s'} &
G_1 \ar[d]^{s} & & G_1\times_{G_0} G_1 \ar[r]^{~~~~~~~~~~~~M(m)} \ar[d]_{t'} &
G_1 \ar[d]^{t}\\
G_1 \ar[r]_{~~~~s}
& G_0 &, & G_1 \ar[r]_{~~~~t}
& G_0 }
\end{displaymath}

 which says that when applying "the composition" $M(m)$ to a pair of arrows $(f,g)$ such that $t(f)=s(g)$, we get an arrow $g\circ f:=M(m)(f,g)$ such that $s(g \circ f)=s(f), t(g \circ f)=t(g)$.

Apply $M$ to the commutative diagram which was used to define $\psi$, we have the commutative diagram
\begin{displaymath}
\xymatrix{
 & & G_1\times_{G_0} G_1\times_{G_0} G_1 \ar[lld]_{p_{12}} \ar@{.>}[d]|-{M(\psi)=(M(m), id)} \ar[drr]^{p_3} & & \\
 G_1\times_{G_0} G_1 \ar[rd]_{M(m)} & & G_1\times_{G_0} G_1 \ar[ld]_{s'} \ar[rd]^{t'} & & G_1 \ar[dl]^{id} \\
 & G_1 \ar[rd]_{t} & & G_1 \ar[dl]^{s} & \\
 & & G_0 & &
 }
\end{displaymath}
where $p_{12}, p_3$ are the obvious projections.
Indeed, $M(\psi)$ is the obvious projection $(M(m),id)$,
since the bottom diagram is cartesian in $\mathfrak{Set}$ and the outer diagram commutes (by the second axiom).
By analogous consideration, we have that $M(\phi)=(id, M(m))$.

Thus, the third diagram yields the commutative diagram
\begin{displaymath}
\xymatrix{
G_1\times_{G_0} G_1\times_{G_0} G_1 \ar[r]^{(M(m),id)} \ar[d]_{(id, M(m))} &
G_1\times_{G_0} G_1 \ar[d]^{M(m)} \\
G_1\times_{G_0} G_1 \ar[r]_{M(m)}
& G_1 }
\end{displaymath}
which expresses the associativity of the composition $M(m)$, i.e., $h\circ(g\circ f)=(h\circ g)\circ f$ for any triple $(f,g,h)$ of arrows with $t(f)=s(g),s(h)=t(g)$.

The 4-th, 5-th diagrams yield the commutative diagrams
\begin{displaymath}
\xymatrix{
G_0  \ar[r]^{M(e)} \ar[rd]_{id} & G_1 \ar[d]^{s} & & G_0  \ar[r]^{M(e)} \ar[rd]_{id} & G_1 \ar[d]^{t}\\
& G_0 &, & & G_0
}
\end{displaymath}
which say that "the unit map" $M(e)$ assigns an arrow $id_a:=M(e)(a) \in G_1$ with $S(id_a)=a=t(id_a) $ to each vertex $a \in G_0$.

Similar arguments for showing $M(\psi)=(M(m),id)$ show that

$$M(\delta)=(M(e),id):G_1 \cong G_0 \times_{G_0} G_1 \rightarrow G_1 \times_{G_0} G_1$$,
$$M(\rho)=(id,M(e)):G_1 \cong G_1 \times_{G_0} G_0 \rightarrow G_1\times_{G_0} G_1$$.

Thus, the last diagram yields the commutative diagram

\begin{displaymath}
\xymatrix{
G_1 \cong G_0\times_{G_0} G_1 \ar[r]^{(M(e), id)} \ar[dr]_{id} &
G_1\times_{G_0} G_1 \ar[d]^m & G_1 \times_{G_0} G_0 \cong G_1 \ar[l]_{(id,M(e))} \ar[dl]^{id}\\
 & G_1 &
 }
\end{displaymath}
which says that $f \circ id_a=f$ for any $(a,f) \in G_0 \times G_1$ with $s(f)=a$ and $g=id_b \circ g$ for any $(g,b) \in G_1 \times G_0$ with $t(g)=b$.

All of these say that $(G,M(m),M(e))$ is a category.

For the converse, given a category $\CC$, define the functor $M:{\LC}_{\mathfrak C} \rightarrow \mathfrak{Set}$ by the following;

$M(G)=\mathfrak{Gr}(G,U(\CC))$ for $G \in ob({\LC}_{\mathfrak C})=ob(\mathfrak{Gr}^{op}_f)$,

$M(\alpha)=\mathfrak{Gr}(\alpha,U(\CC))$ for morphisms $\alpha$ in $\mathfrak{Gr}_f$,

$M(m):U(\CC)_1 \times_{U(\CC)_0} U(\CC)_1 \rightarrow U(\CC)_1, (f,g) \mapsto g \circ f$,

$ M(e):U(\CC)_0 \rightarrow U(\CC)_1, a \mapsto id_a $.

Then, all diagrams commute obviously.
Finally, one can easily check that two constructions are mutually inverse.

\end{proof}

\begin{remark}
\rm For the Lawvere theory ${\LC}_{\mathfrak C}$ we have defined, the functors $U_{{\LC}_{\mathfrak C}}$ and $F_{{\LC}_{\mathfrak C}}$ coincide with forgetful functor and free construction of proposition \ref{free}.
\end{remark}

In a similar way we can show that groupoids are models for a Lawvere ${\mathfrak Gr}$-theory.

\begin{definition}
\label{groupoid}
${\cal L}_{\mathfrak G}$ is the Lawvere ${\mathfrak Gr}$-theory having the following presentation:

generators:  $ m:\overrightarrow{2}\rightarrow\overrightarrow{1}$, $e:\overrightarrow{0} \rightarrow \overrightarrow{1}$, $\iota:\overrightarrow{1}\rightarrow\overrightarrow{1}\}$

 axioms(relations): all those appearing in definition \ref{cat} plus
\begin{displaymath}
\xymatrix{
\vec{1} \ar[r]^{\iota} \ar[dr]_{t^{op}} & \vec{1} \ar[d]^{s^{op}} & & \vec{1} \ar[r]^{\iota} \ar[dr]_{s^{op}} & \vec{1} \ar[d]^{t^{op}}  & & \vec{1} \ar[r]^{\xi} \ar[d]_{t^{op}} & \vec{2} \ar[d]^m & \vec{1} \ar[l]_{\zeta} \ar[d]^{s^{op}} \\
 & \vec{0} &, &  & \vec{0} &, & \vec{0} \ar[r]_e & \vec{1} & \vec{0} \ar[l]^e
 }
\end{displaymath}

where $\xi$ and $\zeta$ are the unique morphisms in ${\LC}_{\mathfrak G}$ making the following diagrams in ${\LC}_{\mathfrak G}$ commute
\begin{displaymath}
\xymatrix{
\vec{1} \ar@/_/[ddr]_{\iota} \ar@/^/[drr]^{id} \ar@{.>}[dr]|-{\xi} & & & & \vec{1} \ar@/_/[ddr]_{id} \ar@/^/[drr]^{\iota} \ar@{.>}[dr]|-{\zeta}  \\
& \vec{2} \ar[d]^{s'^{op}} \ar[r]_{t'^{op}} & \vec{1} \ar[d]^{s^{op}} & & & \vec{2} \ar[d]^{s'^{op}} \ar[r]_{t'^{op}} & \vec{1} \ar[d]^{s^{op}} \\
& \vec{1} \ar[r]_{t^{op}} & \vec{0} &, & & \vec{1} \ar[r]_{t^{op}} & \vec{0}
}
\end{displaymath}
\end{definition}
Note that such unique morphisms $\xi, \zeta$ exist in ${\LC}_{\mathfrak G}$, since the bottom diagrams are cartesian in ${\LC}_{\mathfrak G}$ and the outer diagrams commute.

\begin{theorem}
The category $Mod({\cal L}_{\mathfrak G})$ of ${\cal L}_{\mathfrak G}$-models is equivalent to the category ${\mathfrak Grpd}$ of groupoids.
\end{theorem}
\begin{proof}
Following the proof of theorem \ref{catt}, we have that for any model $M$ there exists a graph $G\in {\mathfrak Gr}$ such that $M\circ J={\mathfrak Gr}(i_,G)$.

We refer to the proof of theorem \ref{catt} for what concerns those diagrams already appearing there.

The first and second diagrams in definition \ref{groupoid} yield the following diagram in ${\mathfrak Set}$
\begin{displaymath}
\xymatrix{
G_1 \ar[r]^{M(\iota)} \ar[dr]_{t^{op}} & G_1 \ar[d]^{s^{op}} & & G_1 \ar[r]^{M(\iota)} \ar[dr]_{s^{op}} & G_1 \ar[d]^{t^{op}}\\
 & G_0 &, &  & G_0
 }
\end{displaymath}
which say that the ``inverse map" $M(\iota)$ assigns to any arrow $f\in G_1$ an arrow $f^{-1}:=M(\iota)(f)\in G_1$ such that $s(f^{-1})=t(f)$ and $t(f^{-1})=s(f)$.

Applying $M$ to the the commutative diagram defining $\xi$ we obtain another commutative diagram
\begin{displaymath}
\xymatrix{
G_1 \ar@/_/[ddr]_{M(\iota)} \ar@/^/[drr]^{id}
\ar@{.>}[dr]|-{M(\xi)=(M(\iota),id)} \\
& G_1\times_{G_0} G_1 \ar[d]^{s'^{op}} \ar[r]_{t'^{op}}
& G_1 \ar[d]^{s^{op}} \\
& G_1 \ar[r]_{t^{op}}
& G_0
}
\end{displaymath}
$M(\xi)$ is $(M(m),id)$,
since the bottom diagram is cartesian in $\mathfrak{Set}$ and the outer diagram commutes (by the second axiom).
By analogous considerations, we have that $M(\zeta)=(id, M(\iota))$.

Therefore the third diagram yields the commutative diagram
\begin{displaymath}
\xymatrix{
G_1 \ar[r]^{(M(\iota),id)} \ar[d]_{t^{op}} & G_1\times_{G_0} G_1 \ar[d]^{M(m)}  & G_1 \ar[l]_{(id, M(\iota))} \ar[d]^{s^{op}}\\
G_0 \ar[r]_{M(e)} & G_1 & G_0 \ar[l]^{M(e)}
}
\end{displaymath}
which says that $f\circ f^{-1}=id_{t(f)}$ and $f^{-1}\circ f=id_{s(f)}$.

These, together with what proved in theorem \ref{catt}, say that $(G,M(m),M(e),M(\iota))$ is a groupoid.

For the converse, as in the proof of theorem \ref{catt}, given a groupoid ${\cal G}$, using the inclusion ${\mathfrak Grpd}\subset{\mathfrak Cat}$ to apply the forgetful functor $U$ to ${\cal G}$, define the functor $M:{\LC}_{\mathfrak G} \rightarrow \mathfrak{Set}$ by the following:

$M(G)=\mathfrak{Gr}(G,U({\cal G}))$ for $G \in ob({\LC}_{\mathfrak G})=ob(\mathfrak{Gr}^{op}_f)$,

$M(\alpha)=\mathfrak{Gr}(\alpha,U({\cal G}))$ for morphisms $\alpha$ in $\mathfrak{Gr}_f$,

$M(m):U({\cal G})_1 \times_{U({\cal G})_0} U({\cal G})_1 \rightarrow U({\cal G})_1, (f,g) \mapsto g \circ f$,

$ M(e):U({\cal G})_0 \rightarrow U({\cal G})_1, a \mapsto id_a $.

$M(\iota):U({\cal G})_1\rightarrow U({\cal G})_1, f \mapsto f^{-1} $.

Then all diagrams commute.
Finally, one can check that two constructions are mutually inverse.
\end{proof}

\section{Finitely presentable categories and models}

We want now to prove that finitely presentable objects are just finitely presentable models for a Lawvere ${\mathfrak Gr}$-theory.

In this section, $\LC$ will denote a Lawvere ${\mathfrak Gr}$-theory where $\mathfrak{Gr}$ is considered as a category, i.e., a $\mathfrak{Set}$-category.
Recall that an object $C$ in a category $\CC$ is \emph{finitely presentable} if the representable functor
$${\CC}(C,-):\CC \rightarrow {\mathfrak Set}$$
preserves filtered colimits.

\begin{definition}
\label{presentation}
A model $M\in Mod(\LC)$ is finitely presentable when there exist $G$ and $H$ in ${\mathfrak Gr}_f$ such that $M$ is the coequalizer
\begin{displaymath}
\xymatrix{
F(H) \ar@<1ex>[r]^\alpha \ar@<-1ex>[r]_\beta & F(G) \ar[r]^q & M
 }
\end{displaymath}
\end{definition}

We call this a finite presentation of $M$.

\begin{proposition}
\label{reflective}
$Mod(\LC)$ is a reflective subcategory of $[{\LC},{\mathfrak Set}]$.
\end{proposition}
\begin{proof}
See \cite{LR}.
\end{proof}

This implies in particular that $Mod(\LC)$ is complete and cocomplete.

\begin{lemma}
\label{representable}
${\cal L}(G,-)=F(iG)$ for $G\in{\mathfrak Gr}_f$.
\end{lemma}
\begin{proof}
Our statement says that for a model $M$
$$Mod({\cal L})({\cal L}(G,-),M)={\mathfrak Gr}(iG,U(M))$$
but this follows from proposition 4.1 of \cite{NP}.
\end{proof}

\begin{proposition}
\label{section}
Free models on finite graphs form a dense family of generators of $Mod(\LC)$.
\end{proposition}
\begin{proof}

By proposition \ref{reflective} $Mod(\LC)$ is a reflective subcategory of $[{\LC},{\mathfrak Set}]$; in $[{\LC},{\mathfrak Set}]$ every model $M$ is the colimit of representable functors ${\cal L}(JG,-)$ for $G$ finite; these, on the other hand, are in $Mod(\LC)$ as, by lemma \ref{representable}, ${\cal L}(G,-)=F(iG)$ for $G\in{\mathfrak Gr}_f$; so the colimit  $M$ exists in $Mod(\LC)$.
\end{proof}

\begin{proposition}
\label{finitedense}
If $M$ is a finitely presentable model, then it admits a presentation (a coequalizer as in definition \ref{presentation}) such that the $q$, as graph morphism, admits a section $s$, that is, $q\circ s=id_M$ in ${\mathfrak Gr}$
\begin{displaymath}
\xymatrix{
F(H) \ar@<1ex>[r]^\alpha \ar@<-1ex>[r]_\beta & F(G) \ar[r]^q & M \ar@/^/[l]^{s}
 }
\end{displaymath}
\end{proposition}
\begin{proof}
Let $M$ be a finitely presentable model and take a presentation of it
\begin{displaymath}
\xymatrix{
F(H') \ar@<1ex>[r]^{~~\alpha'} \ar@<-1ex>[r]_{~~\beta'} & F(G') \ar[r]^{~~~q'} & M
 .}
\end{displaymath}
Consider the following adjunctions of $\alpha', \beta'$
\begin{displaymath}
\xymatrix{
H' \ar@<1ex>[r]^{\alpha''~~} \ar@<-1ex>[r]_{\beta''~~} & UF(G')
. }
\end{displaymath}
Let $R_0$ be the smallest equivalence relation containing $<\alpha'(v),\beta'(v)>$, for $v\in |H'|$, and, since $|UF(G')|=|G'|$, let $r:G'\rightarrow G'/R_0$ be the quotient morphism. Applying $F$ we get a morphism $F(r):F(G')\rightarrow F(G'/R_0)$. Note that $F(r)$ is an epimorphism, because $r$ is and $F$ is left-adjoint to $U$. We can now define a morphism $\bar{q}:F(G'/R_0)\rightarrow M$: it acts on an equivalence class of $F(G'/R_0$ as $q$ acts on a representative, and this is well-defined because of how $R_0$ is defined; it acts on morphisms precisely as $q$ does, as $R_0$ is an equivalence relation just on objects.
\begin{displaymath}
\xymatrix{
 & & F(G'/R_0) \ar[d]^{\bar{q}} \ar[r]^{~~~~~p} & N \\
F(H') \ar@<1ex>[r]^{~~\alpha'} \ar@<-1ex>[r]_{~~\beta'} & F(G') \ar[ru]^{F(r)} \ar[r]_{~~~q'} & M \ar[ru]_t
 }
\end{displaymath}
We have that $\bar{q}\circ F(r)\circ\alpha=q\circ\alpha=q\circ\beta=\bar{q}\circ F(r)\circ\beta$ and we want to show that
\begin{displaymath}
\xymatrix{
F(H) \ar@<1ex>[r]^{F(r)\circ\alpha~~~~~} \ar@<-1ex>[r]_{F(r)\circ\beta~~~~~} & F(G/R_0) \ar[r]^{~~~~~\bar{q}} & M \ar@/^/[l]^{~~~~s}
 }
\end{displaymath}
 is a coequalizer. It remains to prove the universal property. So let $(N,p)$ such that $p\circ F(r)\circ\alpha=p\circ F(r)\circ\beta$. By universality we have a unique morphism $t:M\rightarrow N$ such that $p\circ F(r)=t\circ q'$. Since $q'=\bar{q}\circ F(r)$ we have that $p\circ F(r)=t\circ \bar{q}\circ F(r)$, and, since $F(r)$ is an epimorphism, we get that $p=t\circ\bar{q}$.
 Observe now that, since $F(G/R_0)$ and $M$ are graphs with same vertexes, there exists a section $s:M\rightarrow F(G/R_0)$ to $\bar{q}$.
Note finally that $H'$ is finite by assumption and $F(G/R_0)$ is finite since $G$ is and $R_0$ just identifies some vertexes.
\end{proof}

\begin{proposition}
\label{finpres}
The finitely presentable models form a dense family of generators in $mod({\cal L})$, stable under finite colimits, and every model is a filtered colimit of finitely presentable ones.
\end{proposition}
\begin{proof}
The proof with parallel that proposition 3.8.12 in \cite{B}. Let $\F$ be the full subcategory of finitely presentable models. For a model $M$ consider the overcategory $\F/M$ and the forgetful functor $\phi:\F/M\rightarrow Mod(\LC)$. Following \cite{B} and using proposition \ref{finitedense}, we have ${\rm colimit}\phi=(M,s_{(F,f)})$, where $s_{(F,f)}=f:\phi((F,f))=F\rightarrow M$.

That the colimit above is cofiltered, that is, that $F/M$ is cofiltered, follows from the fact that $F$ is stable in $Mod(\LC)$ under finite colimits. Let us prove this. Following \cite{B}, we soon have that $\F$ is stable under finite coproducts. It is stable also under coequalizers. The proof is again similar to that in \cite{B}, however we need to apply proposition \ref{section}. Suppose $P$ and $Q$ are finitely presentable, let $u,v:P\rightarrow Q$ be two morphism, and let $(R,r)$ be the coequalizer: we want to prove that $R$ is also finitely presentable. Since $P$ and $Q$ are finitely presentable we can consider the diagram
\begin{displaymath}
\xymatrix{
F(H) \ar@<1ex>[d]^b \ar@<-1ex>[d]_a & F(K) \ar@<1ex>[d]^d \ar@<-1ex>[d]_c  & \\
F(G) \ar@<1ex>@{.>}[r]^x \ar@<-1ex>@{.>}[r]_y \ar@{>>}[d]_p & F(J) \ar@{>>}[d]^q &  \\
P \ar@<1ex>[r]^u \ar@<-1ex>[r]_v & Q \ar@{>>}[r]_r \ar@/^/[u]^{s} & R
}
\end{displaymath}
the existence of the lifts $x$ and $y$ of respectively $u$ and $v$ is a consequence of proposition \ref{section}, since we can choose a presentation of $Q$ admitting a section $s:Q\rightarrow F(J)$ of $q$. The proof follows now as in cite{B}, showing that $R$ admits indeed a presentation
\begin{displaymath}
\xymatrix{
F(G\amalg K) \ar@<1ex>[r]^{~~~x\amalg c} \ar@<-1ex>[r]_{~~~y\amalg d} & F(J) \ar[r]^{~r\circ q} & R
 }
\end{displaymath}
\end{proof}

\begin{lemma}
\label{freegen}
Free models on finite graphs are finitely presentable models.
\end{lemma}
\begin{proof}
Let $F(G)$ be a free model with $G$ finite and consider a cofiltered colimit $X={\rm colim}X_i$, then by adjointness
$$Mod({\cal L})(F(G),{\rm colim}X_i)={\mathfrak Gr}(G,U({\rm colim}X_i)$$
since $U$, being finitary monadic (see proposition \ref{monadic}) preserves filtered colimits, we have
$${\mathfrak Gr}(G,U({\rm colim}X_i)={\rm colim}{\mathfrak Gr}(G,U(X_i))$$
finally, since $G$ is finitely presentable
$${\rm colim}{\mathfrak Gr}(G,U(X_i))={\rm colim}Mod({\cal L})(F(G),X_i)$$
thus free finitely presentable models are finitely presentable objects.
\end{proof}

Before enouncing the main result, the following one is expected, having started our construction with finitely presentable categories:

\begin{proposition}
$Mod({\cal L})$ is locally finitely presentable.
\end{proposition}
\begin{proof}
$Mod({\cal L})$ is cocomplete by proposition \ref{reflective}. Free generators are finitely presentable by lemma \ref{freegen} and by proposition \ref{finitedense} form a dense, thus strong, family of generators.
\end{proof}

%
We conclude with the main result:

\begin{theorem}
\label{equivalence}
Finitely presentable models correspond to finitely presentable categories.
\end{theorem}
\begin{proof}
Let $M$ a finitely presentable model and take a presentation
\begin{displaymath}
\xymatrix{
F(H) \ar@<1ex>[r] \ar@<-1ex>[r] & F(G) \ar[r] & M
 }
\end{displaymath}
since $F(H)$ and $F(G)$ are finite presentable objects, and since these are stable under finite colimits,
it follows that $M$ is a finitely presentable object.

For the converse, suppose that for $M\in Mod({\cal L})$ we have an isomorphism
$$Mod({\cal L})(M,{\rm colim}X_i)\cong {\rm colim}Mod({\cal L})(M,X_i)$$
for any filtered colimit $X={\rm colim}X_i$. By proposition \ref{finpres}, $M$ is a filtered colimit of finitely presentable ones: $(M,s_{(F,f)})={\rm colim}\phi(F,f)$; so, substituting, we obtain
$$Mod({\cal L})(M,M)\cong {\rm colim}Mod({\cal L})(M,\phi(F,f))$$
Let $f:M\rightarrow F$ be the morphism corresponding to the identity on $M$: together with $s_{(F,f)}$ expresses $M$ as a retract of $P$ and so $M$ as a coequalizer of $(id_F,f\circ s_{(F,f)}):F\rightarrow F$. By proposition \ref{finpres}, $M$ is finitely presentable.

\end{proof}


\begin{thebibliography}{100}\frenchspacing\small

\bibitem[Bo94]{B} F. Borceux, \emph{Handbook of categorical algebra 2, Categories and structures,} Cambridge University Press, Cambridge (1994).

\bibitem[K82]{K} G.M. Kelly, \emph{Basic concepts of enriched category theory,} \rm London Mathematical Society lecture note series, 64, (1982).

\bibitem[LP11]{LP} S. Lack, J. Power, \emph{Gabriel-Ulmer duality and Lawvere theories enriched over a general base,} \rm Preprint, (2011).

\bibitem[LR11]{LR} S. Lack, J. Rosicky, \emph{Notions of Lawvere theory,} \rm Applied Categorical Structures, 19(1), (2011), 363-391.

\bibitem[L64]{L} W.F. Lawvere, \emph{Functorial Semantics of Algebraic Theories,} \rm PhD Thesis, (1964).

\bibitem[NP09]{NP} K. Nishizawa, J. Power, \emph{Lawvere theories enriched of a general base,} \rm J. Pure Appl. Algebra, 213(3), (2009), 377-386.

\bibitem[P99]{P} J. Power, \emph{Enriched Lawvere theories,} \rm Theory and Applications of Categories, 6(7), (1999), 83-93.

\end{thebibliography}
\end{document}